\documentclass[11pt,a4paper]{article}
\usepackage[utf8]{inputenc}
\usepackage{amsmath}
\usepackage{amsfonts}
\usepackage{amssymb}
\usepackage{float}
\usepackage{url}

\usepackage{epsf,array,delarray,amsmath,amssymb,amsthm}
\usepackage{graphicx}
\usepackage{cite}

\usepackage{color} 

\newcommand{\half}{  {\scriptstyle \frac{1}{2} }\,  }
\newcommand{\R}{   {\mathbb R}  }
\newcommand{\Z}{   {\mathbb Z}  }

\newcommand{\sF}{{\mathcal F}}

\newcommand{\sT}{{\mathcal T}}
\newcommand{\bb}{\beta}

\newcommand{\sgn}{ \hbox{\rm sgn} }

\newcommand{\parens}[1]{\left(#1\right)}

\theoremstyle{plain}
\newtheorem{proposition}{Proposition}
\newtheorem{theorem}{Theorem}

\author{Philip Ernst, L.C.G. Rogers, and Quan Zhou}
\title{The value of foresight}
\begin{document}
\maketitle

\begin{center}
\textit{We dedicate this work to our colleague, mentor, and friend, Professor Larry Shepp (1936-2013)}
\end{center}

\begin{abstract}
Suppose you have one unit of stock, currently worth 1,  which you must sell before time $T$. The Optional Sampling Theorem tells us that whatever stopping time we choose to sell, the expected discounted value we get when we sell will be 1. Suppose however that we are able to see $a$ units of time into the future, and base our stopping rule on that; we should be able to do better than expected value 1. But how much better can we do? And how would we exploit the additional information? The optimal solution to this problem will never be found, but in this paper we establish remarkably close bounds on the value of the problem, and we derive a fairly simple exercise rule that manages to extract most of the value of foresight.
\end{abstract}

\section{Introduction.}\label{intro}
What is the value of foresight in a financial market? This is a question that intrigued Larry Shepp (see page 2 of \cite{SheppSem}) and seems an interesting question in the context of insider trading; if we could know one minute in advance what the price of a stock was going to do, what would we be prepared to pay for that information? Of course, it is rather fanciful to imagine that we could possibly be told the price of the stock at some time in the future, but we might imagine a situation where some market participants received information only after a delay, which would confer the same kind of advantage on those who got the information earlier. In a modern financial market, any such differences would be measured in microseconds, a timescale on which conventional models of stock prices could not be trusted, but one of the first observations of this paper is that the value of foresight can be equivalently interpreted in terms of the value of a fixed-window lookback option; at a time of your choosing, you may sell the stock for the best price which the stock achieved in the previous $a$ units of time. This transforms the question into an American option pricing problem, but not one that is possible to solve in closed form, since the state variable at time $t$ is the entire history of the stock from time $t-a$ to time $t$ -a path-valued state. Moreover, even if we were to discretize time, the state vector will be high dimensional, so existing numerical methods will struggle to cope. Nevertheless, recent developments allow good progress to be made on the question, as we shall see. \\
\indent To begin with, we set some notation. We shall take the sample space to be the path space $C(\R^+,\R)$ with canonical process $W$, Wiener measure $P$ and the usual $P$-augmentation $(\sF_t)$ of the canonical filtration of $W$. We denote by $\sT$ the class of all $(\sF_t)$-stopping times. We fix some $a>0$ which represents the foresight available to the insider. The insider may choose to stop at any stopping time of the larger filtration $(\sF'_t) \equiv (\sF_{t+a})$; we denote the set of all $(\sF'_t)$-stopping times by $\sT'$. 

 At an abstract level, this set-up can be considered an example of {\em grossissement}  ---  the enlargement of filtrations. This theory was developed from the late 1970's on, starting with the works of Barlow \cite{Barlow}, \cite{Bayraktar}, Jeulin \& Yor \cite{Yor1}, \cite{Yor2}, \cite{Yor3}, and further developed  by others including Yoeurp \cite{Yoeurp}, and by It\^{o}'s extension of the stochastic integral (see \cite{Ito}). The topic has since continued to flourish, primarily because of its natural connection with insider trading. Varied formulations of the insider trader's advantage over other agents in a financial market are addressed in \cite{Amendinger}, \cite{Back}, \cite{Elliott}, \cite{Pontier}, \cite{Imkeller3}, \cite{Imkeller2}, \cite{Imkeller1}, \cite{Jacod}, \cite{Karatzas}, \cite{Mansuy}, and \cite{Protter}, among others. It has to be understood that the theory of enlargement of filtrations is not a {\em universal} theory; results are only established for particular classes of enlargement, such as filtrations enlarged with an honest time, or filtrations with an initial enlargement. The results proved say that if the enlargement has one of these particular structural forms, then any $(\sF_t)$-local martingale is a $(\sF'_t)$-semimartingale, and the semimartingale decomposition is then identified.  The following proposition shows that none of these general results can be applied to the problem we consider here.
 
 \begin{proposition}\label{prop0}
 The process $W$ is not a semimartingale in the filtration $(\sF'_t)$.
 \end{proposition}

\begin{proof}
Consider the simple integrands
\begin{equation}
H^n_t \equiv n^{-1/2} \sum_{j=1}^n \sgn(\Delta^n_j)\; I_{  \{ (j-1)a < nt \leq ja \}},
\label{Hn}
\end{equation}
where 
\begin{equation}
\Delta^n_j \equiv W( ja/n ) - W( (j-1)a/n).
\label{Deltadef}
\end{equation}
The processes $H^n$ are left-continuous, bounded, and $(\sF'_t)$-previsible; indeed, $H^n_t$ is measurable on $\sF'_0 \equiv \sF_a$. Now consider the (elementary) stochastic integral
\begin{equation}
H^n \cdot W = n^{-1/2} \sum_{j=1}^n |\Delta^n_j| = n^{-1} \sum_{j=1}^n \sqrt{n}|\Delta^n_j| \;.
\end{equation}
The random variables $\sqrt{n} \Delta^n_j$ are independent zero-mean gaussians with common variance $a$.  By the Weak Law of Large Numbers, $H^n \cdot W$ converges in probability to $E|W_a|$ as $n \rightarrow \infty$. But the Dellacherie-Bichteler theorem (see, for example, Theorem IV.16.4 in \cite{RW2}) says that $W$ is a semimartingale if and only if whenever a sequence $H^n$ of bounded previsible simple processes tends uniformly to zero, then the simple stochastic integrals $H^n\cdot W$ tend to zero in probability. We conclude that $W$ is not an $(\sF'_t)$-semimartingale.

\end{proof}

The message from Proposition \ref{prop0} is that none of the results from enlargement of filtrations will help us here --- the problem addressed is {\em concrete, challenging, and not amenable to general theory} --- which is why it appealed to Larry Shepp.

So what {\em are} we able to do? To begin with, applying the methods of \cite{bobs}, we obtain remarkably tight upper and lower bounds on the value of foresight. This methodology is based on simulations, so there is no simple interpretation of the exercise rule which arises. However, in Section \ref{S3} we develop simple and transparent rules based on heuristic arguments which we are then able to compare with the bounds from Section \ref{S2}; the resulting (explicit) rules essentially achieve the lower bound which comes from the simulation approach of \cite{bobs} applied in Section \ref{S2}.

\subsection*{Preliminaries.}
To start with, we set out some notation and make various standardizations of the question.  There is a fixed time horizon $T>a>0$ by which time the investor must have sold the stock. The stock price process $S$ will be the solution to an SDE
\begin{equation}
dS_t = S_t( \sigma dW_t + \mu_t dt)I_{\{t \leq T\} }, \qquad S_0 = 1,
\label{dS}
\end{equation}
where $\mu$ is a bounded $(\sF_t)$-previsible process, and $\sigma>0$ is a constant. Notice that $S_t = S_T$ for all $t \geq T$, and where necessary we make the convention that $S_t = 1$ for $t<0$.  It is immediate from the definitions that $\tau' \in \sT'$ if and only if $\tau'+a \in \sT$. We let $Q$ denote the pricing measure, given in terms of the Cameron-Martin-Girsanov martingale
\begin{equation}
d\Lambda_t = \Lambda_t (r-\mu_t) \; dW_t, \qquad \Lambda_0 = 1,
\end{equation}
by
\begin{equation}
\frac{dQ}{dP}\biggr\vert_{\sF_t} = \Lambda_t.
\end{equation}

What is the time-0 value of stopping at $\tau' \in \sT'$?  The obvious immediate answer to this question is just $E^Q[ \exp(-r\tau') S_{\tau'}]$, but is it clear that there is no issue arising from the fact that $S_{\tau'}$ is not $\sF_{\tau'}$-measurable?  We can see that this is in fact correct by the following argument. At time $\tau'$, the investor receives $S_{\tau'}$ which he can then place in the bank account for $a$ units of time, so that at $(\sF_t)$-stopping  time $\tau = \tau'+a$ he has $e^{ra} S_{\tau'} = e^{ra} S_{\tau-a}$. This random variable is $\sF_\tau$-measurable, so the time-0 value of it is just given by the usual expression, 
\begin{equation*}
E^Q[ e^{-r\tau} e^{ra} \,S_{\tau-a}] = E^Q[ e^{-r\tau'}\, S_{\tau'}],
\end{equation*}
as expected.  So what we find is that the time-0 value of stopping at $\tau' \in \sT'$ is the $Q$-expectation of the discounted value of the stock at the time of exercise. We may as well therefore work with the discounted stock, and work in the pricing measure $Q$. Equivalently, we may (and shall) assume for simplicity that
\begin{equation}
r = \mu = 0.
\label{simple1}
\end{equation}
Since the choice of $\sigma$ amounts to the choice of a time unit, we may and shall assume that 
\begin{equation}
\sigma = 1.
\label{simple2}
\end{equation}
Our model for the asset dynamics is therefore 
\begin{equation}
S_t = \exp( W_t - \half t) \equiv \exp( X_t),
\label{eq2}
\end{equation}
where $X$ is the drifting Brownian motion
\begin{equation}
X_t = W_t + ct,
\label{Xdef}
\end{equation}
with the special value\footnote{Later on, we shall derive expressions for various probabilities and expectations associated with $X$, and it turns out to be notationally cleaner to work with a general drift, which is why we write $X$ as \eqref{Xdef}.} $c = -\half$. Notice that we have by convention fixed $S_0 = 1$. \\

 What we want to understand then is 
\begin{equation}
v(a) \equiv \sup_{\tau' \in \sT',\; 0 \leq  \tau' \leq T} E[\; S_{\tau'}\; ] \equiv \sup_{\tau \in \sT, \;  0 \leq \tau \leq T+a}
E[\; S_{\tau-a}  \; ].
\label{eq3}
\end{equation}
It is clear that $v$ will be increasing, and $v(0) = 1$, but our aim is to determine as accurately as possible what $v(a)$ is, and to identify a good approximation to the optimal stopping time.  The first item of business, dealt with in Section \ref{S2}, is to show that the value $v$ can be alternatively expressed as
\begin{equation}
v(a) = \sup_{\tau \in \sT, \; 0 \leq \tau \leq T} E[ \; Z_\tau \; ],
\label{eq7}
\end{equation}
where\footnote{... using the convention that $S_u = 1$ for $u<0$, and $S_u = S_T$ for $u \geq T$ ...} 
\begin{equation}
Z_t \equiv \sup\{ S_u: t-a \leq u \leq t \}.
\label{Zdef}
\end{equation}
In other words,  the value $v(a)$ is the value of an {\em American fixed-window lookback option.}

\section{Foresight as lookback.}\label{S2}
Recalling the convention that $S_u = 1$ for $u<0$, and $S_u = S_T$ for $u \geq T$, we have a simple proposition.

\begin{proposition}\label{prop1}
With $\tau$ denoting a generic $(\sF_t)$-stopping time,
\begin{equation}
v(a) \equiv \sup_{a \leq \tau \leq T+a} E[\; S_{\tau-a}  \; ] = \sup_{0 \leq \tau \leq T} E[ \; Z_\tau \; ],
\label{eq9}
\end{equation}
where $Z_t \equiv \sup\{ S_u: t-a \leq u \leq t \}$.
\end{proposition}

\vskip 0.1 in \noindent
{\sc Proof.}  Because $S_u = S_T$ for all $u \geq T$, it is clear that $Z_t \leq Z_{t \wedge T}$.  Therefore, for any stopping time $\tau$ such that $ a \leq \tau \leq T+a$, we have
\begin{equation*}
 S_{\tau-a} \leq Z_\tau \leq Z_{\tau \wedge T}.
\end{equation*}
Therefore 
\begin{equation}
v(a) \equiv \sup_{a \leq \tau \leq T+a} E[\; S_{\tau-a}  \; ] \leq \sup_{0 \leq \tau \leq T} E[ \; Z_\tau \; ].
\label{eq10}
\end{equation}
For the reverse inequality, suppose that $\tau$ is a stopping time, $0 \leq \tau \leq T$, and define a new random time $\tilde \tau$ by
\begin{equation}
\tilde\tau	= \inf\{ u \geq \tau\vee a:  S_{u-a} = Z_\tau\}.
\label{tildetau}
\end{equation}
Clearly $\tau \vee a \leq \tilde\tau \leq \tau + a$. We claim that $\tilde\tau$ is a stopping time, as follows:
\begin{eqnarray*}
\{ \tilde\tau \leq v \} &=& \{ \hbox{\rm for some $u \in [\tau \vee a,v], \; 
S_{u-a} = Z_\tau$ }  \}
\\
&=& \{ \tau\vee a \leq v\} \cap \{
\hbox{\rm for some $u \in [(\tau \vee a)\wedge v ,v], \; 
S_{u-a} = Z_{\tau\wedge v}$ }  \}
\\
& \in & \sF_v  \;   ,
\end{eqnarray*}
since the event
$ \{\exists u \in [(\tau\vee a)\wedge v ,v], \; 
S_{u-a} = Z_{\tau\wedge v}  \}$ is $\sF_v$-measurable, as is 
$(\tau\vee a) \wedge v$.  Now we see that
\begin{equation}
Z_\tau = S_{\tilde\tau - a},
\end{equation}
and therefore
\begin{equation}
E[ \; Z_\tau \; ] = E[ \; S_{\tilde\tau - a}  \; ]
\leq \sup_{a \leq \tau \leq T+a} E[\; S_{\tau-a}  \; ],
\label{prop1_conc}
\end{equation}
since  $ a \leq \tilde\tau \leq T+a$.  Since $0 \leq \tau \leq T$ was any stopping time, we deduce that
\begin{equation}
\sup_{0 \leq \tau \leq T} E[ \; Z_\tau \; ] \leq \sup_{a \leq \tau \leq T+a} E[\; S_{\tau-a}  \; ],
\end{equation}
and the proof is complete.
\vskip 0.1 in 
\hfill $\square$

\bigbreak
The importance of Proposition \ref{prop1} is that it turns the problem of calculating the value of foresight into the calculation of an American fixed-window lookback option. By discretizing the time, we shall instead calculate numerically the value of a Bermudan fixed-window lookback option. Of course, we need to account for the difference between American and Bermudan prices, but this is in essence a solved problem; see Broadie, Glasserman \& Kou \cite{BGK}.

\bigbreak
So we shall take the standardized asset dynamics \eqref{eq2}, fix some time horizon $T>0$ which is subdivided into $N_T$ steps of length $h = T/N_T$, and consider the problem of bounding
\begin{equation}
v_h(a) \equiv \sup_{0 \leq \tau_h \leq T} E[ \; Z^{(h)}_{\tau_h} \; ],
\label{vh}
\end{equation}
where $\tau_h$ is a stopping time taking values in the set $h \Z^+$, and 
\begin{equation}
Z^{(h)}_t	= \max\{ S_{kh}: t-a \leq kh \leq t \}.
\label{Zh}
\end{equation}
Since this discretization is now fixed for the rest of the section, we shall drop the appearance of $h$ in the notation and refer to $Z$ for $Z^{(h)}$, $v$ for $v_h$, $\tau$ for $\tau_h$.  We are now exclusively considering the optimal stopping of  a functional of a discrete-time Markov process. A moment's thought shows that the process $Z$ is not  Markov, but the process
\begin{equation}
x_t	= (S_{t-mh},\;  \ldots,\;  S_t)	\qquad (t \in h \Z^+)
\label{xdef}
\end{equation}
is Markovian (where $mh = a$), and the payoff process $Z$ is simply a function of the Markov process $x$:
\begin{equation}
Z_t	= g(x_t),
\end{equation}
where $g(x) \equiv \max\{x_0, \ldots, x_m \}$ is the largest component of the $(m+1)$-vector $x$.  We use the approach of \cite{bobs}, which is a combination of several techniques developed over the last twenty or so years, and in summary consists of the four steps:
\begin{itemize}
\item pretend that the stopping reward process $Z$ is itself Markovian, and by discretizing $Z$ onto a suitably-chosen finite set of values estimate the transition probabilities of this finite state Markov chain by simulation (this is the approach of Barraquand \& Martineau \cite{BM});
\item solve the optimal stopping problem for this finite state Markov chain by dynamic programming;
\item use the solution to generate a stopping rule whose performance is evaluated by simulation;
\item use the dual characterization of the value of the problem (see \cite{R1}, \cite{HK}, \cite{AB}) to find a hedging martingale.  
\end{itemize}
The method is fully explained in \cite{bobs}, and illustrated with examples, so there is no need to discuss it further here, except to highlight one point, which is used in various places in \cite{bobs} and is needed here.  Suppose that $(M_t)_{t \geq 0}$ is any strictly positive martingale, $M_0=1$. Then we may equivalently express the value 
\begin{eqnarray}
v(a) &\equiv& \sup_{0 \leq \tau \leq T} E[ \; Z_\tau \; ]
\nonumber
\\
&=& \sup_{0 \leq \tau \leq T} E[ \; M_\tau \, (Z_\tau/M_\tau) \; ]
\nonumber
\\
&=& \sup_{0 \leq \tau \leq T} \tilde E[ \;Z_\tau/M_\tau \; ],
\label{numeraire}
\end{eqnarray}
where the probability $\tilde P$ equivalent to $P$ is defined by using the likelihood-ratio martingale $M$:
\begin{equation}
\frac{d\tilde P}{dP}\biggr\vert_{\sF_t} = M_t.
\label{dPtildedP}
\end{equation}
In the present application, it is natural to use $M_t = S_t$ as the change-of-measure martingale, and the effect of this is to change the Brownian motion $W$ into a Brownian motion with drift 1. Thus when we do simulations, we simulate a Brownian motion with drift 1 in place of $W$, and the stopping reward process is changed to $Z_t/S_t$.  This is a good thing to do, because the value of $Z_t$ will be close to $S_t$, so the binning procedure of the Barraquand-Martineau approach should achieve a lot better accuracy than we would get without this measure transformation. 
{Moreover, after the change of measure the stopping rule depends on $(Z_t/S_t, t)$ instead of $(S_t, t)$, which is more sensible and resembles the stopping rule we will propose below.}
We do indeed find that the accuracy is substantially improved by doing this change of measure; the results are reported in Table 1.

 \small
\begin{table}[h!]
\begin{center}
\begin{tabular}{cccccc}
\hline
$a/h$ & Lower  & SE(low) & Upper & SE (up)& Gap(\%)\\
\hline
1 & 1.054 & 0.43 & 1.055 & 0.53&0.09 \\
2 & 1.074 & 0.62 & 1.076 & 0.77&0.18 \\
3 & 1.088 & 0.73 & 1.092 & 0.90&0.37 \\
4 & 1.100 & 0.83 & 1.104 & 1.0&0.37 \\
5 & 1.109 & 0.95 & 1.114 & 1.1&0.45 \\
6 & 1.117 & 1.0 & 1.123 & 1.2&0.53 \\
7 & 1.123 & 1.1 & 1.131 & 1.2&0.71 \\
8 & 1.129 & 1.2 & 1.137 & 1.3&0.70 \\
9 & 1.135 & 1.3 & 1.144 & 1.4&0.79 \\
10 & 1.140 & 1.4 & 1.149 & 1.4&0.78 \\
11 & 1.144 & 1.4 & 1.154 & 1.5&0.87 \\
12 & 1.149 & 1.5 & 1.159 & 1.6&0.86 \\
13 & 1.152 & 1.6 & 1.164 & 1.7&1.03 \\
14 & 1.156 & 1.6 & 1.168 & 1.7&1.03 \\
15 & 1.159 & 1.7 & 1.172 & 1.8&1.11 \\
16 & 1.163 & 1.8 & 1.175 & 1.8&1.02 \\
17 & 1.166 & 1.8 & 1.179 & 1.9&1.10 \\
18 & 1.168 & 1.9 & 1.182 & 1.8&1.18 \\
19 & 1.171 & 1.9 & 1.185 & 1.9&1.18 \\
20 & 1.174 & 2.0 & 1.188 & 1.9&1.18 \\
\hline
\end{tabular}
{\caption{Upper and lower bounds of $\sup_{\tau< T} E[Z_\tau]$ from simulation using method of \cite{bobs} with $h = 1/2500$, $N_T = 250$. `SE' is the standard error, in basis points (1e-4). Simulation parameters: number of bins = 200, number of samples per bin = 200, lower-bound simulation = 50,000, upper-bound simulation = 10,000, sub-simulation per step = 50.}}\label{table:bermuda}
\end{center}
\end{table}
\normalsize

\section{Explicit stopping rules.}\label{S3}
As was stated earlier, the stopping rules which are derived in Section \ref{S2} are the output of a simulation; they have no particular structure or interpretation, and a different simulation run will generate a different stopping rule.  The methodology of \cite{bobs} is generic, and works just the same for essentially any Markov process, and any stopping function of the process, but we may hope to improve in individual applications by exploiting the specific structure of that application, which is what we shall do here. 

\bigbreak
We work in terms of the log price $X_t = W_t + ct$, and use the notations
\begin{equation}
\bar X_t \equiv \sup_{0 \leq u \leq t} X_u, \qquad
\bar X_{[t-a,t]} \equiv \sup_{t-a \leq u \leq t} X_u \equiv \log Z_t\; .
\label{notation}
\end{equation}
The first time we may need to consider stopping is the stopping time
\begin{equation}
 \tau_0 \equiv \inf \{\;  t: \bar X_{[t-a,t]} = X_{t-a} \; \},
 \label{tau0def}
\end{equation}
because up until that time we have $X_{t-a} < \bar X_{[t-a,t]}$, so it will be suboptimal to exercise - waiting a little longer may improve, and will not make the reward less.  But at time $\tau_0$, continuation means we have to let go of the good value $X_{\tau_0-a}$ in the hope of doing better in the future, and we may in fact do worse. Whether we should optimally continue will depend on the entire path of $X$ from $\tau_0-a$ to $\tau_0$, but we will simplify to consider only rules where continuation is decided by the value of $X_{\tau_0} - X_{\tau_0 -a}\;$; if this is higher than some threshold $q <0$ we shall continue, otherwise we stop.

An important observation is the fact that
\begin{equation}
\bar X_{[\tau_0-a,\tau_0]}  = \bar X_{\tau_0},
\label{simplify}
\end{equation}
as a moment's thought will reveal. Therefore we have
\begin{equation}
\tau_0 = \inf \{\;  t: \bar X_t = X_{t-a} \; \}.
\label{tau0def_bis}
\end{equation}

Now it is clear that the choice of threshold value $q$ will depend on the time to go when we have to make the stop/continue decision, and this makes complete solution of the problem much more complicated. So what we shall do is to propose a {\em modified} problem, where the time by which we must sell is not a fixed time $T$ but an independent random time $\alpha \sim \exp(\eta)$ for some $\eta>0$. This gives us a renewal property that allows us to make progress, and obtain explicit expressions. When it comes to converting the solution to the modified problem into an exercise rule for the original problem, what we do is set the threshold according to the value of $\eta$ which makes the expectation of $\alpha$ equal to the time to go.

For this modified problem, we define a stopping rule $R(q)$ depending on the chosen threshold as follows:
\begin{itemize}
\item[(i)] Wait until $\tau_0 \wedge \alpha$;
\item[(ii)] If $\alpha < \tau_0$, stop and receive $Z_\alpha = \exp(\bar X_\alpha)$;
\item[(iii.a)] If $\tau_0 < \alpha$ and $X_{\tau_0} - \bar X_{\tau_0} < q$, stop and receive $Z_{\tau_0}= \exp(\bar X_{\tau_0} )$;
\item[(iii.b)] If $\tau_0 < \alpha$ and $X_{\tau_0} - \bar X_{\tau_0} > q$, forget and continue.
\end{itemize}
In the final eventuality (iii.b), `forget and continue' means that we wipe away the whole path of $X$ in the time interval $[0, \tau_0)$, keeping only the value $X_{\tau_0}$, and  restart the rule from that point.  If the value of this strategy is $K$, then we have the identity
\begin{eqnarray}
K &=& E[ \;\exp(\bar X_\alpha) : \alpha < \tau_0]
+ E[\;  \exp(\bar X_{\tau_0}) : \tau_0 < \alpha, X_{\tau_0} - \bar X_{\tau_0} < q \;  ]
\nonumber
\\
&& \qquad\qquad\qquad  + K E[ \; \exp(X_{\tau_0}): \tau_0 < \alpha, X_{\tau_0} - \bar X_{\tau_0} > q  \;  ].
\label{Keq}
\end{eqnarray}
It is therefore apparent that we can evaluate this particular stopping rule provided we can find explicit expressions for
\begin{eqnarray}
A_0 &=&  E[ \;\exp(\bar X_\alpha) : \alpha < \tau_0\;],  \label{A0}
\\
A_-(q) &=& E[\;  \exp(\bar X_{\tau_0}) : \tau_0 < \alpha, X_{\tau_0} - \bar X_{\tau_0} < q \;  ],   \label{A-}
\\
A_+(q) &=& E[ \; \exp(X_{\tau_0}): \tau_0 < \alpha, X_{\tau_0} - \bar X_{\tau_0} > q  \;  ].  \label{A+}
\end{eqnarray}
We can. If $\varphi$ denotes the density of the standard $N(0,1)$ distribution, and $\Phi$ denotes its distribution function, we have the following result.

\begin{proposition}\label{prop2}
For $q<0$, denote
\begin{eqnarray}
\bb	&=& \sqrt{c^2+2\eta},
\label{bdef}
\\
\nu &=& \frac{2}{\sqrt{a}}\; \varphi(\bb\sqrt a) + 2\bb\Phi(\bb\sqrt a) - c - \bb,
\label{nudef}  \\
\nu_a &=&  \frac{2}{\sqrt{a}} \; \varphi(c\sqrt{a}) - 2c\, \Phi(-c\sqrt a),
\label{nuadef}  \\
\nu_\alpha &=&  2\bb(\Phi(\bb\sqrt a) - \half) -2c(\Phi(c\sqrt a) - \half)
 \nonumber\\
 &&\qquad\qquad \qquad\quad 
 + \frac{2}{\sqrt a} \bigl\lbrace\; \varphi(\bb\sqrt a) - \varphi(c\sqrt a)
 \;\bigr\rbrace ,
 \label{nualphadef}
 \nonumber \\
\psi_0(q)&=& e^{-\eta a} \biggl\lbrace\;
\frac{2}{\sqrt a} \; \varphi( (q-ca)/\sqrt a) - 2c\, \Phi( (q-ca)/\sqrt a)
\;\biggr\rbrace,
\label{psi0def}
\\
\psi_1(q) 
&=& e^{(c+\half-\eta)a} \biggl[ \; \frac{2}{\sqrt a}\;\bigl\lbrace
 \; \varphi(\bar c \sqrt a)
-\varphi( (q-\bar c a)/\sqrt a) \;\bigr\rbrace
\nonumber \\
&& \qquad\qquad\qquad + 2\bar c\bigl\lbrace \; \Phi(\bar c \sqrt a) -
 \Phi( -(q-\bar c a)/\sqrt a)  \; \bigr\rbrace \; \biggr],
 \label{psi1def}
\end{eqnarray}
where $\bar c \equiv 1 + c$.  Then:
\begin{eqnarray}
A_0 &=& \frac{\nu_\alpha+ (1-e^{-\eta a}) \nu_a}{\nu-1},
\label{A0_ans}
\\
A_-(q) &=& \frac{\psi_0(q)}{\nu -1},
\label{A-ans}
\\
A_+(q) &=& \frac{\psi_1(q)}{\nu -1}.
\label{A+ans}
\end{eqnarray}
\end{proposition}

\noindent
{\sc Proof.} The process $Y_t \equiv X_t - \bar X_t$ is a diffusion taking values in the negative half-line, and reflecting from zero. The process $\bar X_t$ is its local time at zero, and the path of $Y$ can be decomposed into a Poisson process of excursions, as It\^o \cite{itopp} explained. For a more extensive discussion of excursion theory, in particular, the notion of marked excursion processes, see \cite{GTTE} or Chapter VI.8 of \cite{RW2}.

The process $X$ runs until $\tau_0 \wedge \alpha$, which happens in the first excursion of $Y$ which {\em either} lasts at least $a$, {\em or} contains an $\eta$-mark.  Let $n$ denote the excursion measure, $\zeta$ denote the lifetime of an excursion, $P^c$ denote the law of $X_t = W_t +ct$, and $H_b$ denote the first time $X$ hits $b$.
{
The excursion law can be characterized as the limit as $\varepsilon \downarrow 0$ of the (rescaled) law of $X$ started at $ -\varepsilon <0$ run until it hits zero; see VI.50.20 in \cite{RW2}. The rescaling required is to multiply by 
$n (f: \inf f < -\epsilon )$, and it is known that $ n (f: \inf f < -\epsilon )\sim \epsilon^{-1}$ (see VI.51.2 in \cite{RW2}). We therefore have}
\begin{eqnarray}
n(\zeta>a) &=& \lim_{\varepsilon\downarrow 0}\varepsilon^{-1}P^c[\; H_0 >a \;|\; X_0 = -\varepsilon]  \nonumber
\\
&=&\lim_{\varepsilon\downarrow 0}\varepsilon^{-1}P^c[\; \bar X_a < \varepsilon \;|\; X_0 = 0]  \nonumber
\\
&=& \lim_{\varepsilon\downarrow 0}\varepsilon^{-1}
\biggl\lbrace\;
\Phi\biggl( \frac{\varepsilon-ca}{\sqrt{a}}  \biggr)
- e^{2c\varepsilon} \Phi\biggl( \frac{-\varepsilon-ca}{\sqrt{a}}  \biggr)
\;\biggr\rbrace
\label{eq29}
\\
&=& \frac{2}{\sqrt{a}} \; \varphi(c\sqrt{a}) - 2c\, \Phi(-c\sqrt a),
\label{eq30}
\\
&\equiv & \nu_a,  \nonumber
\end{eqnarray}
where \eqref{eq29} is a standard result on Brownian motion; see, for example, \cite{BS}. By similar reasoning, 
\begin{eqnarray}
n(\alpha < \zeta < a) &=&\lim_{\varepsilon\downarrow 0}\varepsilon^{-1}
\int_0^a P^c[ H_0\in ds | X_0 = -\varepsilon] \; (1-e^{-\eta s})
\nonumber\\
&=& \lim_{\varepsilon\downarrow 0}\varepsilon^{-1}
\int_0^a \frac{\varepsilon e^{-\varepsilon^2/2s}}{\sqrt{2\pi s^3}}\;
e^{c\varepsilon- c^2 s/2} (1-e^{-\eta s})  \; ds
\nonumber\\
&=& \int_0^a \frac{ e^{-c^2 s/2}}{\sqrt{2\pi s^3}}\;
 (1-e^{-\eta s})  \; ds
 \nonumber\\
 &=&  2\bb(\Phi(\bb\sqrt a) - \half) -2c(\Phi(c\sqrt a) - \half)
 \nonumber\\
 &&\qquad\qquad \qquad\quad + \frac{2}{\sqrt a} \bigl\lbrace\; \varphi(\bb\sqrt a) - \varphi(c\sqrt a)
 \;\bigr\rbrace  \label{eq31}
 \\
 &\equiv & \nu_\alpha  \nonumber
\end{eqnarray}
where we recall that $\bb \equiv \sqrt{c^2 + 2 \eta}$.
The rate of excursions which {\em either} last at least $a$, {\em or} contain an $\eta$-mark is 
\begin{equation}
n(\zeta >a, \; \hbox{\rm or} \; \alpha < \zeta)
= n(\alpha < \zeta<a) + n(\zeta >a)
\label{eq33}
\end{equation}
and this is simply the sum of the two expressions \eqref{eq30} and \eqref{eq31}, which is therefore known explicitly. In fact, a few calculations confirm that it is the expression $\nu$ defined at \eqref{nudef}.  Immediately from It\^o excursion theory:
\begin{itemize}
\item $\bar X_{\tau_0 \wedge \alpha} \sim \exp( \;
\nu
\;)$;
\item $
P[ \; \alpha < \tau_0 \;]
= \{ \;\nu_\alpha+ (1-e^{-\eta a}) \nu_a\; \}/\nu
$;
\item $P[\; \tau_0 < \alpha \;] = \nu_a\, e^{-\eta a} /\nu$;
\end{itemize}
and from this the expression \eqref{A0_ans} for $A_0$ follows.

To deal with $A_\pm$, we need to find the measure of excursions which get to time $a$ without killing, and which are  in $dy$  at time $a$. We use the reflection principle and the Cameron-Martin-Girsanov theorem to derive
\begin{eqnarray}
g(y)dy &\equiv& n(a < \zeta\wedge\alpha, Y_a \in dy) 
\label{gdef}
\\
&=&  \lim_{\varepsilon\downarrow 0} \varepsilon^{-1} P^c[ \;
H_0\wedge \alpha >a, Y_a \in dy \; \vert \; Y_0 = -\varepsilon
\;]    \nonumber
\\
&=&\lim_{\varepsilon\downarrow 0} \varepsilon^{-1} e^{c(y+\varepsilon)-c^2a/2
-\eta a}
\{ \;   \varphi( (y + \varepsilon)/\sqrt{a}) - \varphi( (y  -\varepsilon)/\sqrt{a}) 
\;\} dy /\sqrt{a} 
\nonumber
\\
&=&  \frac{-2y}{a^{3/2}} \; \varphi((y-ca)/\sqrt a)\; e^{-\eta a} \; dy.
\label{g_ans}
\end{eqnarray}
Straightforward calculations lead us to 
\begin{eqnarray*}
\int_{-\infty}^q g(y) \; dy
&=& e^{-\eta a} \biggl\lbrace\;
\frac{2}{\sqrt a} \; \varphi( (q-ca)/\sqrt a) - 2c\, \Phi( (q-ca)/\sqrt a)
\;\biggr\rbrace,
\\
&\equiv & \psi_0(q),  \\
\int_q^0 e^y g(y) \; dy
&=& e^{(c+\half-\eta)a} \biggl[ \; \frac{2}{\sqrt a}\;\bigl\lbrace
 \; \varphi(\bar c \sqrt a)
-\varphi( (q-\bar c a)/\sqrt a) \;\bigr\rbrace
\nonumber \\
&& \qquad\qquad\qquad  2\bar c\bigl\lbrace \; \Phi(\bar c \sqrt a) -
 \Phi( -(q-\bar c a)/\sqrt a)  \; \bigr\rbrace \; \biggr],
\\
& \equiv & \psi_1(q).
\end{eqnarray*}
Hence $P( Y_{\tau_0} < q \;\vert \;\tau_0 < \alpha ) = \psi_0(q)/\psi_0(0)$. We therefore have
\begin{eqnarray*}
A_-(q) &=& E[\;  \exp(\bar X_{\tau_0}) : \tau_0 < \alpha, X_{\tau_0} - \bar X_{\tau_0} < q \;  ]
\\
&=& \frac{\nu}{\nu-1} \; P[\;\tau_0 < \alpha, X_{\tau_0} - \bar X_{\tau_0} < q \;  ]
\\
&=& \frac{\nu}{\nu-1}\cdot \frac{\nu_ae^{-\eta a}}{\nu}\cdot \frac{\psi_0(q)}{\psi_0(0)}
\\
&=& \frac{\psi_0(q)}{\nu -1}
\end{eqnarray*}
after some simplifications. This is the form of $A_-$ claimed at \eqref{A-ans}.

Finally, we deal with $A_+(q)$. Observe that
\begin{equation}
E[\; \exp( \, Y_{\tau_0} \, )I_{\{ Y_{\tau_0}>q\}}\; \vert \; \tau_0 < \alpha\;]
=\frac{\psi_1(q)}{n(a < \zeta\wedge\alpha)} \; .
\end{equation}
Hence we have 
\begin{eqnarray*}
A_+(q) &=& E[ \; \exp(X_{\tau_0}): \tau_0 < \alpha, X_{\tau_0} - \bar X_{\tau_0} > q  \;  ]
\\
&=& E[\;  \exp(\bar X_{\tau_0}) \;] \cdot
E[ \; \exp(Y_{\tau_0}): \tau_0 < \alpha, Y_{\tau_0} > q  \;  ]
\\
&=& \frac{\nu}{\nu-1}\cdot  P(\tau_0 < \alpha ) \cdot \frac{\psi_1(q)}{n(a < \zeta\wedge\alpha)}
\\
&=& \frac{\psi_1(q)}{\nu -1}
\end{eqnarray*}
after some simplifications. This is the form claimed at \eqref{A+ans}.

\vskip 0.1 in \hfill $\square$

Now we can draw all the pieces together. Proposition \ref{prop2} was a step on the way to evaluating $K$, the value of the proposed strategy; we found $K$ expressed in terms of $A_0$, $A_-(q)$ and $A_+(q)$ at \eqref{Keq}, and it remains just to write the answer cleanly.

\begin{theorem}\label{thm1}
The value of the rule $R(q)$ is given by
\begin{equation}
K = \frac{A_0 + A_-(q)}{1-A_+(q)} = \frac{\nu_\alpha +(1-e^{-\eta a}) \nu_a + \psi_0(q)}{\nu-1-\psi_1(q)}  \;  .
\label{K_ans}
\end{equation}
The denominator $\nu - 1 - \psi_1(q)$ is always positive.
\end{theorem}

\medskip\noindent
{\sc Proof.} The expression \eqref{K_ans} is a trivial arrangement of \eqref{Keq}, so all that remains is to deal with the final assertion.  Using the fact that $c = -\half$, we see that 
\begin{equation*}
\psi_1(-\infty) = \;
{\uparrow\lim}_{q \downarrow -\infty} \psi_1(q)  
= e^{-\eta a}\biggl\lbrace \; \frac{2}{\sqrt a}\; \varphi(\sqrt{a}/2) -\Phi(-\sqrt{a}/2)\;
\biggr\rbrace.
\end{equation*}
Therefore
\begin{eqnarray*}
\nu - 1- \psi_1(q) &>& \nu - 1 - \psi_1(-\infty)
\\
&=& \nu_a + \nu_\alpha - 1 - \psi_1(-\infty)
\\
&>& \nu_a - 1 - \psi_1(-\infty)
\\
&=& \frac{2}{\sqrt a}\; \varphi(\sqrt{a}/2) + \Phi(\sqrt{a}/2)- 1 - \psi_1(-\infty)
\\
&=& (1-e^{-\eta a}) \biggl\lbrace \; \frac{2}{\sqrt a}\; \varphi(\sqrt{a}/2) -\Phi(-\sqrt{a}/2)\;
\biggr\rbrace
\\
&>& 0.
\end{eqnarray*}

\vskip 0.1 in \hfill $\square$

\subsection*{Using the stopping rule $R(q)$.}
\begin{figure}[H]
\begin{center}
\includegraphics[width=0.49\linewidth]{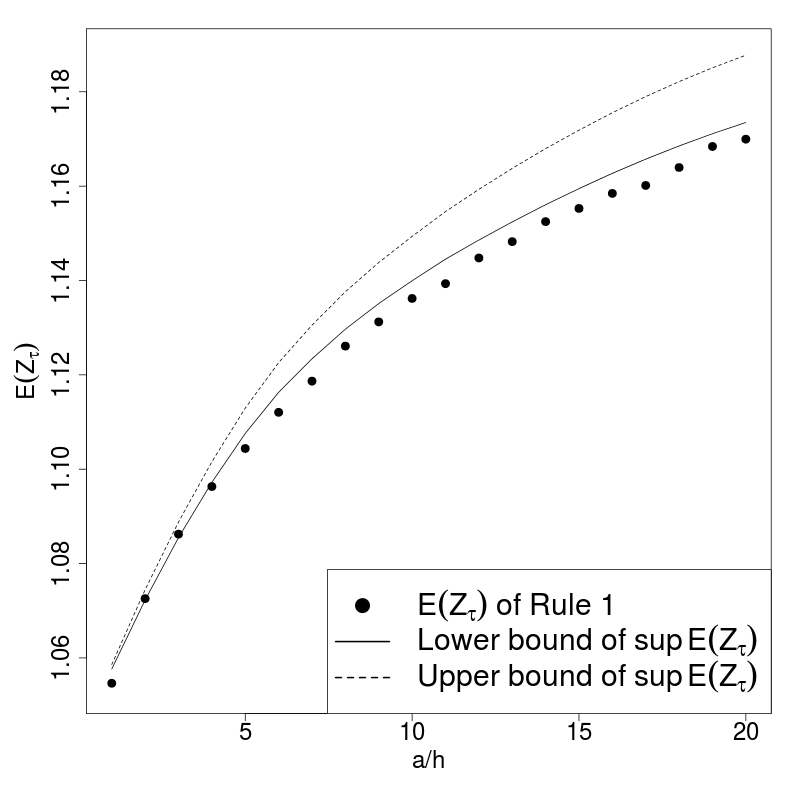}
\includegraphics[width=0.49\linewidth]{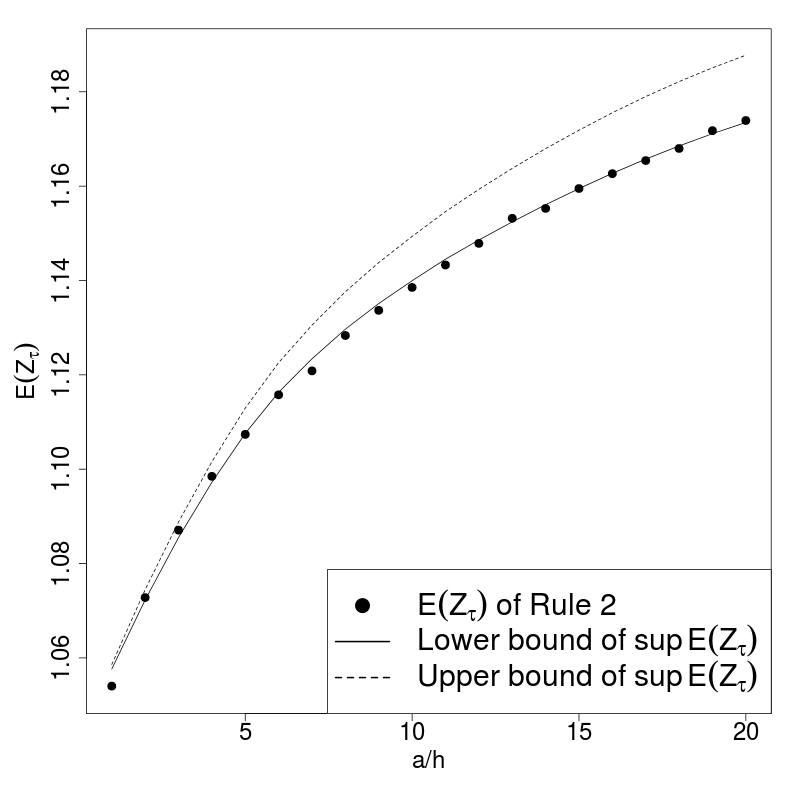} 
\caption{The simulation results with $h = 1/2500$ and $N_T = 250$.}\label{fig:sim}
\end{center}
\end{figure}

{\bf Rule 1.}
How is the stopping rule $R(q)$ analysed in the preceding section relevant to the original problem? Holding $\eta$ fixed, there will be an optimal $q^* = q^*(\eta)$ which maximizes the value $K$ given by \eqref{K_ans}. 
Now recall how the stopping rule works; we let the process run until the stopping time
\begin{equation*}
\tau_0 = \inf \{\;  t: \bar X_{[t-a,t]} = X_{t-a} \; \} = 
 \inf \{\;  t: \bar X_t = X_{t-a} \; \},
\end{equation*}
and at that moment we stop if and only if $Y_{\tau_0} \equiv X_{\tau_0}
-\bar X_{\tau_0} < q$, else we forget and continue.  What value of $q$ do we use? A natural choice is to take
\begin{equation}
q = q^*(\;  (T- \tau_0)^{-1}\;).
\label{q_choice}
\end{equation}
This is because at time $\tau_0$ there is time $(T-\tau_0)$ still to go, and an exponential random variable with rate $(T-\tau_0)^{-1}$ has this as its mean. The result of using this rule is shown in the left-hand panel of Figure \ref{fig:sim}. The dots (evaluated by simulating 50,000 runs of the rule) are visibly close to (but below) the lower bound we obtained by the Barraquand-Martineau technique of \cite{bobs}.
\\

\noindent{\bf Rule 2.}
Can we do better than this? Indeed we can. Firstly, we observe that
\begin{equation}
K(q^*(\eta)) = \exp( -q^*(\eta)).
\label{Kq*}
\end{equation}
This is because when we arrive at time $\tau_0$, we have to choose between stopping and receiving $S_{\tau_0-a}$, or continuing and receiving $K(q^*(\eta))\cdot S_{\tau_0}$ in expectation. Optimal behaviour requires us to continue if and only if
\begin{equation*}
K(q^*(\eta)) \cdot S_{\tau_0} > S_{\tau_0-a}.
\end{equation*}
On the other hand, the rule $R(q^*(\eta))$ requires us to continue if and only if
\begin{equation*}
S_{\tau_0} > S_{\tau_0-a} \cdot \exp( q^*(\eta)).
\end{equation*}
Since $q^*(\eta)$ was chosen optimally, \eqref{Kq*} must therefore hold.
\\

Secondly, if $a = T$, then the value will just be the expectation of the overall maximum, $E \exp( \bar X_T )$.  Routine but tedious calculations  give 
\begin{equation}
E[\; \exp( \bar X_a) \;] = \lambda(a) \equiv \parens{2 + \dfrac{a}{2}}\Phi\parens{\dfrac{\sqrt{a}}{2}} + \sqrt{a} \varphi \parens{\dfrac{a}{4}}. \label{smax.distr}
\end{equation}
If we had arrived at time $\tau_0$ and it turns out that $\tau_0 = T-a$, then by  forgetting and continuing we will actually receive expected reward $\lambda(a) S_{\tau_0}$, whereas if we used Rule 1 we would think we were going to receive $K(q^*(a^{-1}))\cdot S_{\tau_0} = e^{-q^*(a^{-1})} \, S_{\tau_0}$. This suggests that we modify Rule 1, replacing \eqref{q_choice} by
\begin{equation}
 q = q^*(\;(T-\tau_0)^{-1}\;) - q^*(a^{-1}) - \log\lambda(a).
 \label{q_choice_2}
\end{equation}
By making this modification, the stop/continue decision we make in the event that $T- \tau_0 = a$ will be exactly correct.

Of course, the argument just presented is only a rough heuristic, but if we look at the right-hand panel in Figure \ref{fig:sim} we see the results of using  Rule 2. The dots are now essentially coincident with the lower bounds, which is very encouraging, and argues for the use of Rule 2 rather than Rule 1. This rule is something which can be clearly motivated and precisely specified, in contrast to the randomly-generated rules which come from the Barraquand-Martineau technique in Section \ref{S2}.  One could continue to search for other explicit rules which do even better, but such a study is beyond the scope of this paper.

\begin{table}[H]
\begin{center}
\begin{tabular}{cccc}
\hline
$a/h$ & Rule 1 &  Rule 2 & Lower bound\\
\hline
1 & 1.055 &  1.054& 1.054 \\
2 & 1.073 &  1.073 & 1.074   \\
3 & 1.086 &  1.087 &1.088  \\
4 & 1.096  & 1.098 & 1.100 \\
5 & 1.104  & 1.107  & 1.109\\
6 & 1.112 & 1.116 & 1.117  \\
7 & 1.119 & 1.121 &1.123 \\
8 & 1.126  & 1.128 &1.129  \\
9 & 1.131  & 1.134 & 1.135 \\
10 & 1.136 & 1.139 &1.140 \\
11 & 1.139 & 1.143& 1.144 \\
12 & 1.145  & 1.148 &1.149\\
13 & 1.148  & 1.153 &1.152 \\
14 & 1.152 & 1.155  & 1.156 \\
15 & 1.155  & 1.160 & 1.159 \\
16 & 1.158  & 1.163 & 1.163\\
17 & 1.160  & 1.165 & 1.166 \\
18 & 1.164  & 1.168 & 1.168 \\
19 & 1.168  & 1.172  & 1.171 \\
20 & 1.170 & 1.174 & 1.174 \\
\hline
\end{tabular}
{
\caption{$E[Z_\tau]$ estimates from simulation of Rules 1 and 2 with $h = 1/2500$ and $N_T = 250$, averaged over $50,000$ sample paths. Standard errors are 0.001 (to one significant figure) in all cases.
  }\label{table:rq}
}
\end{center}
\end{table}

\section{Supplemental Materials}
The code for all simulations can be found on \url{http://stat.wharton.upenn.edu/~ernstp/blb.cpp}. The pseudocode is available on  \url{http://stat.wharton.upenn.edu/~ernstp/Bermuda_code.pdf}

\pagebreak
\bibliography{VF}
\bibliographystyle{plain}

\end{document}